\documentclass[10pt]{amsart}
\usepackage{graphicx}%
\usepackage{multirow}%
\usepackage{amsmath,amssymb,amsfonts}%
\usepackage{amsthm}
\usepackage{mathrsfs}%
\usepackage[title]{appendix}%
\usepackage{xcolor}%
\usepackage{textcomp}%
\usepackage{manyfoot}%
\usepackage{booktabs}%
\usepackage{listings}
\usepackage{amssymb}
\usepackage{amsfonts}
\usepackage{amsmath}
\usepackage{epsfig}
\usepackage{indentfirst}
\usepackage[usenames,dvipsnames]{pstricks}
\usepackage{pst-grad}
\usepackage{pst-plot}
\usepackage[margin=1.0in]{geometry}
\usepackage{graphicx}
\usepackage{tikz}
\usepackage{tikz-cd}
\usepackage{float}
\usepackage{tabularx}
\usepackage{mathtools}
\usepackage{bm}
\usepackage{ytableau}
\ytableausetup{boxsize=9pt}

\theoremstyle{thmstyleone}%
\newtheorem{theorem}{Theorem}[section]
\newtheorem{proposition}[theorem]{Proposition}%

\theoremstyle{thmstyletwo}%
\newtheorem{example}{Example}%

\theoremstyle{thmstylethree}%

\newcommand{\C}{\mathbb{C}}

\usetikzlibrary{%
  matrix,%
  calc,%
  arrows%
}

\begin{document}

\title[Hook fusion procedure for direct product of symmetric groups]{Hook fusion procedure for direct product of symmetric groups}
\author{Dimpi KM}\address{Indian Institute of Science Education and Research Thiruvananthapuram, Email: dtyagi20@iisertvm.ac.in}
\author{Geetha Thangavelu}\address{Indian Institute of Science Education and Research Thiruvananthapuram, Email:  tgeetha@iisertvm.ac.in }

\keywords{Yang-Baxter equation, Fusion procedure, Direct product of symmetric groups, Young tableau, Jucy-Murphy elements}
\subjclass{16G10, Secondary 20C30, 05E10, 16D40}

\begin{abstract}
In this work, we derive a new expression for the diagonal matrix elements of irreducible representations of the direct product group $S_r\times S_s$ using Grime's hook fusion procedure for symmetric groups, which simplifies the fusion procedure by reducing the number of auxiliary parameters needed. By extending this approach to the product group setting, we provide a method for constructing a complete set of orthogonal primitive idempotents.
\end{abstract}

\maketitle
\section{Introduction}
The Yang Baxter equation, initially appeard in statistical physics, has played a fundamental role in various fields, including quantum groups, knot theory, integrable systems, and quantum mechanics. A major development from this equation is the fusion procedure, which offers a way to create new solutions from existing ones.

For a finite dimensional complex vector space V, a \textit{fused solution} of the Yang-Baxter equation is built from a given solution $R(u)$ using two sequences of complex numbers $\bm{c} = (c_1,\dots , c_r)$ and $\bm{c'}=(c'_1,\dots ,c'_s)$. The resulting fused solution, denoted $R_{\bm{c,c'}}(u)$, operates on $V^{\otimes n}\otimes V^{\otimes n'}$ and preserves the subspace $W_{\textbf{c}}\otimes W_{\textbf{\underline{c}}}$ of $V^{\otimes n}\otimes V^{\otimes n'}$, where 
\begin{equation}
W_{\bm{c}} := (\prod_{1\leq i< j\leq n}^{\rightarrow} R_{ij}(c_i -c_j))V^{\otimes n}.
\end{equation}

The associated \textit{fusion function}, given by 
\begin{equation}
F(\bm{c}) := \prod_{1\leq i< j\leq n}^{\rightarrow} R_{ij}(c_i -c_j),
\end{equation}
may present singularities, meaning that for certain values of $\bm{c}$, the operator $F(\bm{c})$ may not be well-defined.

Jucys \cite{Ju66} demonstrated that the primitive idempotents of the symmetric group $S_n$ can be derived through a limiting process applied to a rational function. This concept, now widely known as the \textit{fusion procedure}, has been extended to several algebras, such as Hecke algebras, Brauer algebras, and BMW algebras, see \cite{IM08}, \cite{IM10}, \cite{IM14}. Molev introduced an alternative approach to this procedure for $S_n$, which relies on the existence of a maximal commutative subablgebra generated by Jucys-Murphy elements \cite{Mo08}. His method constructs primitive idempotents by iterative evaluation of a specific rational function, and this framework has since been adapted to several groups and algebras.

The fusion procedure plays an important role in representation theory, helping to construct and analyze algebraic structures. The symmetric group $S_n$ acts on $V^{\otimes n}$ by permuting tensor factors, and the fusion function corresponds to a rational function in $\C[S_n]$, helps in the construction of orthogonal primitive idempotents. By selecting appropriate values for $c$, explicit expressions for the diagonal matrix elements of irreducible representations can be obtained, as explored in the work of Jucys \cite{Ju66}, Cherednik \cite{Che86}, and Nazarov \cite{Naz97}, \cite{Naz98}, \cite{Naz04}.

A new version of the fusion procedure, known as the \textit{hook fusion procedure}, was introduced by Grime \cite{Gri05}. This method optimizes the traditional fusion approach by reducing the number of auxiliary parameters required, making it more computationally efficient while maintaining the structure of irreducible representations. Grime's hook fusion technique has been successfully applied to symmetric groups and Hecke algebra of symmetric groups \cite{Gri07}, thus simplifying the calculations. This paper extends the hook fusion procedure from symmetric groups to the direct product $S_r\times S_s $.

\section{Representation theory of $S_r\times S_s$}

Let us fix some notations and recall the representation theory of the direct product of  symmetric groups $S_r\times S_s$.

Consider the group algebra $\C[S_r\times S_s]$. The standard generators of the group $S_r \times S_s$ are $s_1,\dots, s_{r-1}$, $s_{r+1},\dots,$ and $s_{r+s}$. The subalgebra generated by the elements $s_i$, $1\leq i\leq r-1$ is isomorphic to the group algebra $\C[S_r]$ while the subalgebra generated by $s_i$, $r+1\leq i\leq r+s$ is isomorphic $\C[S_s]$. This naturally gives rise to the following chain of subalgebras:
\begin{equation}
C := \C(S_1 \times S_0) \subset \C(S_2 \times S_0) \subset \dots \subset \C(S_r \times S_0)\subset\C(S_r \times S_1)\subset \dots \subset \C(S_r \times S_s) \label{C}.
\end{equation}
This sequence of subalgebras helps in understanding the representation theory of $S_r\times S_s$ step by step.

A \textit{partition} of an integer $n$, denoted as $\lambda$, is a sequence of positive integers $\lambda:=(\lambda_1,\lambda_2,\dots, \lambda_k)$ satisfying $\lambda_1\geq\lambda_2\geq\dots\geq\lambda_k$ and $\lambda_1 + \lambda_2 + \dots +\lambda_k=n$. In this case we write $\lambda\vdash n$ to indicate that $\lambda$ is a partition of $n$. This concept extends naturally to a \textit{bi-partition}, which is a pair of partitions $\bm{\lambda} = (\lambda_1, \lambda_2)$ where each $\lambda_i$ is partition of $r$ and $s$, respectively. Let $\Lambda$ denote the set of all bi-partitions of $n$. We then define $$\Lambda_{r,s} = \{\bm{\lambda}=(\lambda_1, \lambda_2)\in\Lambda \mid \lambda_1 \vdash r, \lambda_2 \vdash s \}.$$

Given a partition $\lambda=(\lambda_1, \dots, \lambda_k)\vdash n$, the corresponding \textit{Young diagram} consists of $k$ rows of boxes, where the $j$-th row contains $\lambda_j$ left-justified boxes. A \textit{Young tableau} is formed by filling these boxes bijectively with the numbers $1$ to $n$. A \textit{standard Young tableau} is a Young diagram in which the entries increase from left to right in each row and from top to bottom in each column.  Each box, called a \textit{node} $\alpha$, is identified by its  coordinates $(i, j)$, and its \textit{content} is defined as: $c(\alpha) = j-i$. In the Young diagram of shape $\lambda$, the $\alpha$-\textit{hook} of a node $\alpha = (i,j)$ is the set of boxes in $\lambda$ such that $\{(i,j'):j'\geq j\}\cup\{(i',j):i'\geq i\}.$ 

For a bi-partition $\bm{\lambda}=(\lambda_1, \lambda_2)\vdash n$ , we associate two Young diagrams, one for each partition $\lambda_1$ and $\lambda_2$. This pair of diagrams is collectively refered  to as the \textit{Young diagram} of the bi-partition. A \textit{bi-tableau}, denoted as $\bm{T}$,  is created by filling both diagrams with numbers $1$ to $n$ in a bijective manner. A \textit{node} in this setting is represented as $(\alpha, k)$, where $\alpha$ is a node in the Young diagram of $\lambda_k$, and $k\in {1,2}$ indicates which tableau it belongs to. The \textit{content} remains $c(\alpha) = j-i$. We set $c_i(\bm{T})$ as the content of the node containing the integer $i$ in the tableau $\bm{T}$ for the discussion.    

A\textit{ standard $\bm{\lambda}$-tableau} $\bm{T}$ is constructed through a sequence of bi-partitions $(\lambda_1 \to \lambda_2 \to \dots \to \lambda_{r+s})$, starting with $\lambda_1 = ((1), \phi)$, and ending at $\lambda_{r+s} =\lambda= (\lambda_r, \lambda_s)$, where each step $\lambda_t$ is obtained by adding a single box. If $t\leq r$, a box is added to the first tableau, and if $t> r$, a box is added to the second tableau.

Alternatively, a tableau $\bm{T}=(T_1, T_2)$ of shape $\bm{\lambda}=(\lambda_1, \lambda_2)$ is called \textit{standard} if $T_1$ contains numbers from $1$ to $r$ and $T_2$ contains numbers $r+1$ to $r+s $, both following the increasing order rules of standard Young tableaux. 

\begin{example}
Consider the partition $\bm{\lambda} = ((2,1),(2))$. A standard $\bm{\lambda}$-tableau $T$ can be constructed through the following sequence of Young diagrams:
$$\bm{T} = 
\left(\scalebox{0.5}{\ydiagram{1}}, \varnothing\right) 
\longrightarrow \left(\scalebox{0.5}{\ydiagram{2}}, \varnothing\right) 
\longrightarrow \left(\scalebox{0.5}{\ydiagram{2,1}}, \varnothing\right) 
\longrightarrow \left(\scalebox{0.5}{\ydiagram{2,1}}, \scalebox{0.5}{\ydiagram{1}}\right)
\longrightarrow \left(\scalebox{0.5}{\ydiagram{2,1}}, \scalebox{0.5}{\ydiagram{2}}\right).
$$

Equivalently, this path corresponds to the following standard $\bm{\lambda}$-tableau $\bm{T}$, $$\bm{T}=\left(\, \begin{ytableau} 1& 2\\3 \end{ytableau}, \begin{ytableau} 4& 5 \end{ytableau}\right).$$
\end{example}

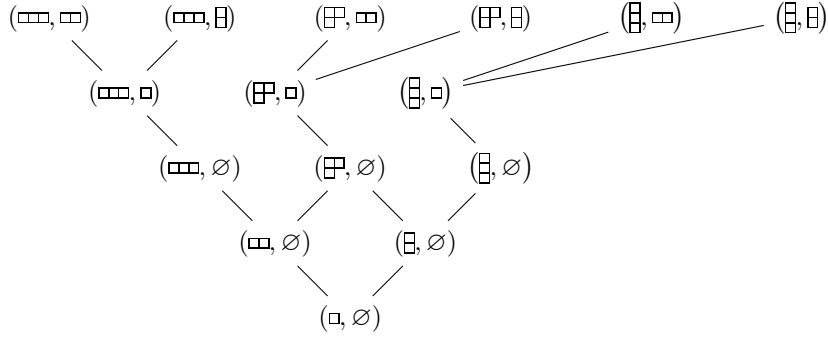
\begin{figure}[h]
    \centering
    \begin{tikzpicture}[baseline={(current bounding box.center)}]

    \node (A1) at (0,0) {$\left(\vcenter{\hbox{\scalebox{0.4}{\ydiagram{3}}}}, \vcenter{\hbox{\scalebox{0.4}{\ydiagram{2}}}}\right)$};
    \node (B1) at (2,0) {$\left(\vcenter{\hbox{\scalebox{0.4}{\ydiagram{3}}}}, \vcenter{\hbox{\scalebox{0.4}{\ydiagram{1,1}}}}\right)$};
    \node (C1) at (4,0) {$\left(\vcenter{\hbox{\scalebox{0.4}{\ydiagram{2,1}}}}, \vcenter{\hbox{\scalebox{0.4}{\ydiagram{2}}}}\right)$};
    \node (D1) at (6,0) {$\left(\vcenter{\hbox{\scalebox{0.4}{\ydiagram{2,1}}}}, \vcenter{\hbox{\scalebox{0.4}{\ydiagram{1,1}}}}\right)$};
    \node (E1) at (8,0) {$\left(\vcenter{\hbox{\scalebox{0.4}{\ydiagram{1,1,1}}}}, \vcenter{\hbox{\scalebox{0.4}{\ydiagram{2}}}}\right)$};
    \node (F1) at (10,0) {$\left(\vcenter{\hbox{\scalebox{0.4}{\ydiagram{1,1,1}}}}, \vcenter{\hbox{\scalebox{0.4}{\ydiagram{1,1}}}}\right)$};

    \node (A2) at (1,-1) {$\left(\vcenter{\hbox{\scalebox{0.4}{\ydiagram{3}}}}, \vcenter{\hbox{\scalebox{0.4}{\ydiagram{1}}}}\right)$};
    \node (B2) at (3,-1) {$\left(\vcenter{\hbox{\scalebox{0.4}{\ydiagram{2,1}}}}, \vcenter{\hbox{\scalebox{0.4}{\ydiagram{1}}}}\right)$};
    \node (C2) at (5,-1) {$\left(\vcenter{\hbox{\scalebox{0.4}{\ydiagram{1,1,1}}}}, \vcenter{\hbox{\scalebox{0.4}{\ydiagram{1}}}}\right)$};

    \node (A3) at (2,-2) {$\left(\vcenter{\hbox{\scalebox{0.4}{\ydiagram{3}}}}, \varnothing\right)$};
    \node (B3) at (4,-2) {$\left(\vcenter{\hbox{\scalebox{0.4}{\ydiagram{2,1}}}}, \varnothing\right)$};
    \node (C3) at (6,-2) {$\left(\vcenter{\hbox{\scalebox{0.4}{\ydiagram{1,1,1}}}}, \varnothing\right)$};

    \node (A4) at (3,-3) {$\left(\vcenter{\hbox{\scalebox{0.4}{\ydiagram{2}}}}, \varnothing\right)$};
    \node (B4) at (5,-3) {$\left(\vcenter{\hbox{\scalebox{0.4}{\ydiagram{1,1}}}}, \varnothing\right)$};

    \node (A5) at (4,-4) {$\left(\vcenter{\hbox{\scalebox{0.4}{\ydiagram{1}}}}, \varnothing\right)$};

    \draw (A1) -- (A2);
    \draw (B1) -- (A2);
    \draw (C1) -- (B2);
    \draw (D1) -- (B2);
    \draw (E1) -- (C2);
    \draw (F1) -- (C2);

    \draw (A2) -- (A3);
    \draw (B2) -- (B3);
    \draw (C2) -- (C3);

    \draw (A3) -- (A4);
    \draw (B3) -- (A4);
    \draw (B3) -- (B4);
    \draw (C3) -- (B4);

    \draw (A4) -- (A5);
    \draw (B4) -- (A5);
    
\end{tikzpicture}
    \caption{Bratteli diagram for $\C[S_3\times S_2]$}
    \label{fig:bratteli}
\end{figure}

The irreducible representations of $S_r \times S_s$ are indexed by elements of $\Lambda_{r,s}$. For each $\bm{\lambda}\in\Lambda_{r,s}$, we denote the corresponding irreducible representation by $V_{\bm\lambda}$. The branching rule for the inclusion chain (\ref{C}) are simple. To illustrate this, a Bratteli diagram can be used, and as an example we consider the case of $\C[S_3\times S_2]$, see Figure \ref{fig:bratteli}.

We define $M(\bm{\lambda})$ to be the set of all standard $\bm{\lambda}$-tableaux. If the Bratteli diagram has multiplicity free edges then, each $\bm{T}\in M(\bm{\lambda})$ corresponds to an one-dimensional subspace $V_{\bm{T}}$ of $V_{\bm{\lambda}}$. The elements of the seminormal basis for $V_{\bm{\lambda}}$ are indexed by standard $\bm{\lambda}$-tableaux. Given a standard $\bm{\lambda}$-tableau $\bm{T}$, let $v_{\bm{T}}$ represents the associated basis vector in $V_{\bm{\lambda}}$. The content of the box in $\bm{T}$ where $i$ is placed is denoted by $c_i(\bm{T})$. The action of the standard generator $s_i$ on $v_{\bm{T}}$ is given by:

\begin{align}
s_i v_{\bm{T}} =
\begin{cases}
 \frac{1}{d_i(\bm{T})} v_{\bm{T}} + \sqrt{1 - d_{i}(\bm{T})^{2}} v_{\bm{T'}} & ; \bm{T'} \text{ is standard} \\
\frac{1}{d_i(\bm{T})} v_{\bm{T}} & ; \bm{T'} \text{ is not standard}
\end{cases}
\end{align}
where $d_i(\bm{T}) = (c_{i+1}(\bm{T})-c_i(\bm{T}))^{-1}$, and $\bm{T'}=s_{i}\bm{T}$ or equivalently $\bm{T'}$ is obtained by exchanging the terms $i$ and $i+1$ in the tableau $\bm{T}.$

The group algebra $\C[S_r\times S_s]$ decomposes into direct sum of matrix algebras, given by
\begin{align}
\C[S_r \times S_s] \cong \bigoplus_{\bm{\lambda} \in \Lambda_{r,s}} \text{End}(V_{\bm{\lambda}}) = \bigoplus_{\bm{\lambda} \in \Lambda_{r,s}} \text{Mat}_{d_{\bm{\lambda}}} (\C),
\end{align}
where $d_{\bm{\lambda}} = \text{dim}(V_{\bm{\lambda}})$. The matrix units $E_{\bm{TT'}}$ in Mat$_{d_{\bm{\lambda}}}(\C)$ are indexed by pairs of standard $\bm{\lambda}$-tableaux $\bm{T}$ and $\bm{T'}$. An explicit isomorphism between these algebras is given by 
\begin{align}
 E_{\bm{TT'}} = \frac{d_{\bm{\lambda}}}{|S_r \times S_s|} F_{\bm{TT'}},
\end{align}
where $F_{\bm{TT'}}$ represents the matrix element associated with the vectors $v_{\bm{T}}$ and $v_{\bm{T'}}$ in the irreducible representation $V_{\bm{\lambda}}$, and is defined as 
\begin{align}
F_{\bm{TT'}} = \sum_{g \in S_r \times S_s} \langle v_{\bm{T}}, g v_{\bm{T'}} \rangle g \in \C (S_r \times S_s).
\end{align}
Now we focus only on diagonal matrix elements, which we denote simply as $E_{\bm{T}} = E_{\bm{TT}}$ and $F_{\bm{T}} = F_{\bm{TT}}$.

We now describe the Jucys-Murphy elements in the group algebra $S_r\times S_s$. These elements are defined as follows:
\begin{align*}
X_1 = 0, \quad X_k = \sum_{i=1}^{k-1} (i\,k), \quad 2 \leq k \leq r,\\
X_{r+1} = 0, \quad X_l = \sum_{i=r+1}^{l-1} (i\,l), \quad r+2 \leq l \leq r+s.
\end{align*}
These elements generate the maximal commutative subalgebra of the group algebra $\C[S_r\times S_s]$. It can be verified that for any $j\in\{1,\dots r, r+1, \dots, r+s\}$, the element $X_j$ commutes with all the elements of the subalgebra positioned at the $(j-1)$th level in the chain (\ref{C}). This property plays a crucial role in understanding the action of Jucys-Murphy elements on the Young basis. 

\begin{proposition}
The Jucys-Murphy element $X_j$ acts on the Young basis element $v_{\bm{T}}$ as 
$$X_j = c_j(\bm{T})v_{\bm{T}},$$
where $\bm{T}$ is the standard $\bm{\lambda}$-tableau.
\end{proposition} 

\begin{proof}
Consider a Young basis element $v_{\bm{T}}\in V_{\bm{\lambda}}$. Since, $V_{\bm{\lambda}} = V_{\lambda_1}\otimes V_{\lambda_2}$, where $V_{\lambda_1}$ and $V_{\lambda_2}$ are the irreducible representations of $\C[S_r]$ and $\C[S_s]$ associated with the partitions $\lambda_1$ and $\lambda_2$, respectively. We can write $v_{\bm{T}}\in V_{\lambda_1}\otimes V_{\lambda_2} $ with $\bm{T}=(T_1, T_2)$. The sum $X_1 + \dots +X_r$ belongs to $ \C[S_r]\otimes 1 \subset \C[S_r\times S_s]$, and corresponds to a symmetric polynomial in the Jucys-Murphy elements of $\C[S_r]$. Thus $$(X_1 + \dots +X_r)v_{\bm{T}} = \sum_{i=1}^{r}c_{i}(T_1)v_{\bm{T}}.$$
This follows from the fact that any symmetric polynomial in the Jucys-Murphy elements of symmetric group acts by scalar multiplication on an irreducible module $V_{\lambda_1}$, with the scalar determined by evaluating the polynomial at the contents of the tableau $T_1$ of shape $\lambda_1$ (See \cite{Mu83}, \cite{DG89}). Similarly, we have $$(X_{r+1} + \dots +X_{r+s})v_{\bm{T}} = \sum_{i=r+1}^{r+s}c_{i}(T_2)v_{\bm{T}}.$$ By definition of the content, we have $c_{\bm{T}}(i)= c_{T_1}(i)$ if $i$ lies in $T_1$, and $c_{\bm{T}}(i)= c_{T_2}(i)$ if $i$ lies in $T_2$. This leads to $$(X_{1} + \dots +X_{j})v_{\bm{T}} = \sum_{i=1}^{j}c_{i}(\bm{T})v_{\bm{T}}.$$ Thus, we conclude that $X_j = c_j(\bm{T})v_{\bm{T}}$, completing the proof.    
\end{proof}

\section{Hook fusion procedure}
We define the subspace $\mathcal{H}_{\bm{T}}$ of $\C^n$ as $$\mathcal{H}_T = \{(z_1, \dots, z_n) \mid z_i = z_j \text{ whenever $i$ and $j$ in the same \textit{principal hook} for $\bm{T}$} \}.$$
Here, the \textit{principal hook} refers to the hook associated with a diagonal node.

Next, we introduce a rational function in the variables $z_1, z_2,\dots, z_n$ that takes values in the group algebra $\C[S_r\times S_s]$ as follows:
\begin{align}
\Phi_{\bm{T}} (z_1, \dots, z_n) :=
\left( \prod_{1 \leq i < j \leq r} \frac{1 - (i j)}{z_i + c_i - (z_j + c_j)} \right)
\left( \prod_{r+1 \leq i < j \leq r+s} \frac{1 - (i j)}{z_i + c_i - (z_j + c_j)} \right),
\end{align}
where we set $c_i:=c_i(\bm{T})$, and $n=r+s$

\begin{theorem}
When restricted to the subspace $\mathcal{H}_{\bm{T}}$, of the rational function $\Phi_{\bm{T}} (z_1, \dots, z_n)$ remains well defined at $z_1 = \dots = z_n$. Furthermore, the value of this restriction at $z_1 = \dots = z_n $ coincides with the diagonal matrix element $F_{\bm{T}}$.
\end{theorem}

\begin{proof}
Since 
\begin{align*}
\Phi_{\bm{T}} (z_1, \dots, z_n) = &\left( \prod_{1 \leq i < j \leq r}(1-\frac{(i\,j)}{z_i + c_i - (z_j + c_j)}) \right)\left( \prod_{r+1 \leq i < j \leq r+s} (1-\frac{(i\, j)}{z_i + c_i - (z_j + c_j)}) \right),
\end{align*}
we can express it as the product
\begin{align*}
\Phi_{\bm{T}} (z_1, \dots, z_n) = \Phi_{T_1}^{1}(z_1,\dots,z_r)\Phi_{T_2}^{2}(z_1,\dots,z_s),
\end{align*}
where $\Phi_{T_{1}}^{1}$ and $\Phi_{T_{2}}^{2}$ are the rational functions in $\C[S_r]$ and $\C[S_s]$, respectively.\\
From the hook fusion procedure for symmetric group algebra \cite{Gri05}, we know that the restriction of $\Phi_{T_1}^{1}(z_1,\dots,z_r)$ to $\mathcal{H}_{T_1}$ and the restriction of $\Phi_{T_2}^{2}(z_{r+1},\dots,z_{r+s})$ to $ \mathcal{H}_{T_2}$ remain regular at $z_1 = z_2 = \dots = z_{r+s}$. Consequently, the restriction of $\Phi_{\bm{T}}(z_1,\dots,z_n)$ to $\mathcal{H}_{\bm{T}}$ is also regular at $z_1=z_2=\dots =z_n$.\\
Now, the diagonal matrix element $ F_{\bm{T}}$ in $\C[S_r\times S_s]$ corresponding to the tableau $T$, can be written as follows: 
$$ F_{\bm{T}} = \left( \sum_{g \in S_r} \langle v_{\bm{T}}, g v_{\bm{T}} \rangle g \right) \left( \sum_{g \in S_s} \langle v_{\bm{T}}, g v_{\bm{T}} \rangle g \right).$$\\ 
This follows from the fact that the action of $s_i$ for $1\leq i\leq r-1$, and $s_k$ for $r+1\leq r+s$ on the Young basis $v_{\bm{T}}$ satisfies $$\langle v_{\bm{T}}, s_ks_i v_{\bm{T}} \rangle = d_kd_i.$$\\
Using the hook fusion procedure for symmetric group, we obtain
$$\Phi_{T_1}^{1} = \sum_{g \in S_n} \langle v_T, g v_T \rangle g, \quad \Phi_{T_2}^{2} = \sum_{g \in S_s} \langle v_T, g v_T \rangle g.$$ This directly implies that $\Phi_{\bm{T}}$ coincides with the diagonal matrix element $F_{\bm{T}}$ in $\C[S_r\times S_s]$.

\end{proof}

\noindent \textbf{Acknowledgements}

The first author's research is supported by Indian Institute of Science Education and Research Thiruvananthapuram PhD fellowship. The second author's research was partially supported by IISER-Thiruvananthapuram, SERB-Power Grant SPG/2021/004200.

\end{document}